\documentclass[12pt]{article}
\usepackage{amsmath,amsthm}
\usepackage{enumitem}
 \usepackage{mathptmx}
\usepackage{graphicx}
\title{Minimal regular graphs with every edge in a triangle}
\author{James Preen\\
\small Cape Breton University, \\[-0.8ex]
\small Sydney, Nova Scotia, \\[-0.8ex] 
\small Canada.\\
\small\tt james\_preen@capebretonu.ca
}
\begin{document}

\newtheorem{lem}{Lemma}[section]
\newtheorem{cor}{Corollary}[section]
\newtheorem{thm}{Theorem}[section]
\newcommand{\dtri}{{\rm diamond}}
\newtheorem{algorithm}{Algorithm}[section]

%\url{https://orcid.org/0000-0002-4688-6930}

\maketitle

\begin{abstract}
Considering regular graphs with every edge in a triangle we prove lower bounds for the number of triangles in such graphs. For $r$-regular graphs with $r\leq 5$ we exhibit families of graphs with exactly that number of triangles and then classify all such graphs using line graphs and even-cycle decompositions. Examples of ways to create such $r$-regular graphs with $r\geq 6$ are also given. In the 5-regular case, these minimal graphs are proven to be the only regular graphs with every edge in a triangle that cannot have an edge removed and still have every edge in a triangle.
\end{abstract}

\section{Introduction}\label{s:intro}

In this paper a triangle in a graph will be defined as 
a set of three distinct vertices $u$, $v$ and $w$ together with three edges $uv$, $uw$ and $vw$. We are interested in graphs with the feature that every edge is in at least one triangle, we refer to this as {\em the triangle property}. 
In particular, we are looking for graphs that have the fewest possible triangles in them with regards their regularity and number of vertices while still having the triangle property.
If $T$ is a triangle including $v$ then we say that $T$ is {\em incident} with $v$.

For 2-regular and 3-regular graphs with the triangle property the only graphs are disjoint copies of complete graphs, $x K_3$ and $z K_4$, respectively, for any positive integers $x$ and $z$. Any disconnected graph with the triangle property must have all components with the property, so henceforth we can suppose all graphs under consideration are connected. 

\begin{thm}\label{t:mintri}
Suppose $G$ is an $r$-regular graph with  
 the triangle property. Every $v\in V(G)$ is incident with at least $\frac{r}{2}$ triangles.
\end{thm}

\begin{proof}
Suppose there are $t$ triangles at a vertex $v$ of degree $r$ in a graph $G$ with the triangle property.
Each of the $r$ edges at $v$ is in a triangle and a triangle incident with $v$ requires two edges from $v$ to neighbours of $v$; note that some of the edges may also be used in other triangles. Thus the $t$ triangles 
require at most $2t$ edges at $v$, and we can conclude that $r\leq 2t$, or $t\geq \frac{r}{2}$, as required.
\end{proof}

An $r$-regular graph with $n$ vertices, the triangle property and $\lfloor\frac{r+1}{2}\rfloor\times \frac{n}{3}$ triangles will be called {\em triminimal}.
Note that for some values of $r$ and $n$ a graph with these exact parameters is impossible, either because of divisibility conditions, or if $n$ is too small with respect to $r$ as we will see in Section \ref{s:lgraph}. 

\begin{cor}\label{c:trimin}
A $r$-regular triminimal graph with $n$ vertices has the fewest triangles amongst $r$-regular graphs  with $n$ vertices and the triangle property.
\end{cor}

\begin{proof}
Let $G$ be any $r$-regular graph with the triangle property.
By Theorem \ref{t:mintri} every vertex $v \in V(G)$ is incident with at least $\frac{r}{2}$ triangles and, when we count the triangles at every vertex, each triangle is counted three times. The number of triangles in $G$ is therefore at least $\frac{r}{2} \times \frac{n}{3}=\frac{rn}{6}$, and note that if $r$ is even then $\lfloor\frac{r+1}{2}\rfloor = \frac{r}{2}$.
If $r$ is odd then $\frac{r}{2}$ is not an integer and so, from Theorem \ref{t:mintri}, we can say there are at least 
$\lfloor\frac{r+1}{2}\rfloor = \frac{r+1}{2}$ triangles at every vertex, giving at least $\frac{(r+1)n}{6}$ triangles for the whole graph.  
\end{proof}

\subsection{Multiple edges}

By our definition of a triangle, multiple edges are permitted in graphs with the triangle property but loops cannot be in a triangle since they must repeat a vertex. Any triangle containing a multiple edge will lead to multiple triangles using those same three vertices, and this will prevent triminimality:

\begin{thm}\label{t:simple}
Any triminimal graph is simple.
\end{thm}

\begin{proof}
From Corollary \ref{c:trimin} all inequalities within the proof must be equalities. For even $r$ every vertex must be incident with exactly $\frac{r}{2}$ triangles and any multiple edge at a vertex would force more triangles.

When $r$ is odd, there must be exactly one edge in two triangles at every vertex to create the necessary $\frac{r+1}{2}$ triangles by Corollary \ref{c:trimin}.
One of the edges at a vertex $v$ could possibly be a double edge instead of an edge in two triangles incident with $v$; for example, we can try to create a triminimal quintic graph with double edges at each vertex by adding extra edges from a 1-factor of a triminimal quartic graph to that graph as is done for the line graph of the cube in Figure \ref{f:multopt}.

\begin{figure}[h]
\begin{center}
\includegraphics[width=0.5\textwidth]{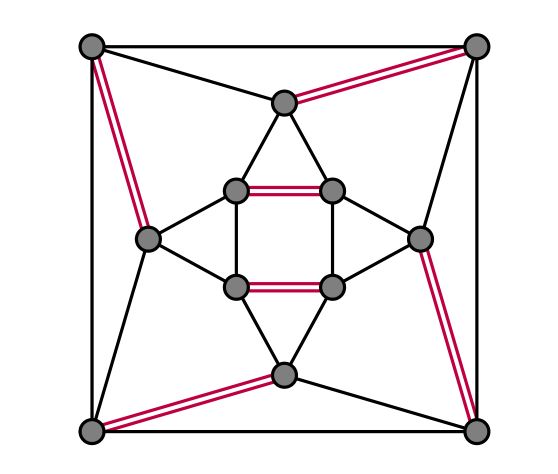}
     \end{center}
\caption{A quintic multigraph with 12 vertices and 14 triangles}\label{f:multopt}  
\end{figure}

In general, for $s\geq 2$, given a $2s$-regular triminimal graph $G$ with $n$ vertices and $\frac{ns}{3}$ triangles we can add $\frac{n}{2}$ edges from a 1-factor of $G$ to give a $(2s+1)$-regular multigraph $H$. 
However, $H$ will have $\frac{ns}{3} + \frac{n}{2} = \frac{n(2s+3)}{6}$ triangles
and the number of triangles in a triminimal graph with $r=2s+1$ is, by definition,
$\frac{n((2s+1)+1)}{6}= \frac{ns}{3}+\frac{n}{3} < \frac{ns}{3}+\frac{n}{2} $. Thus, this construction cannot lead to a triminimal graph.

Henceforth we can assume that in a $(2s+1)$-regular multigraph $G$ there is at least one vertex which is not in a double edge.
Let $w$ be a vertex which is not incident with a double edge but in a triangle $T$ with a double edge. At $w$ there will be more than $s+1$ triangles, contradicting triminimality: two triangles using the vertices of $T$ and at least $\frac{(2s-1)+1}{2}=s$ more from the $2s-1$ neighbours of $w$ not in $T$. 
\end{proof}

Henceforth we will assume all graphs are connected and simple.

\section{Even regularity}\label{s:even}

Recall that, to form the line graph $L(H)$ of a graph $H$, we create a new graph with $|E(H)|$ vertices (corresponding to the edges of $H$) and join pairs of these new vertices 
if the corresponding edges of $H$ had a vertex in common. In particular, if $H$ is a 3-regular graph then $L(H)$ is a 4-regular graph with every vertex incident with at least two triangles.

\begin{cor}\label{c:opt4}
Triminimal 4-regular graphs are line graphs of 3-regular triangle-free graphs.
\end{cor}

\begin{proof}
If $J$ is a triminimal 4-regular graph 
then all inequalities in Theorem \ref{t:mintri} must actually have been equalities for every vertex; there are exactly two triangles incident with every vertex. We can recognise this structure as the line graph $L(H)$ of a 3-regular graph $H$.
 
Line graphs of 3-regular graphs were one of the base families of graphs with the triangle property in \cite{ar:PR}. However, to ensure triminimality, we need to ensure that $H$ is, additionally, triangle-free, as otherwise any triangle in $H$ will also appear in $J$ as a triangle too, and will not come from the three edges incident at a vertex in $H$.
\end{proof}

The line graph construction can be thought of in the following alternate way; we replace each edge of a 3-regular graph $H$ by a path of two edges to form a bipartite graph $B(H)$. The vertices in one part of $B(H)$ have degree 3 and those in the other part have degree 2; such a bipartite graph is called (3,2)-biregular in \cite{ar:bibireg}. We can then form $L(H)$ from $B(H)$ by deleting all of the vertices of degree 3 after adding edges between their neighbours to make triangles; this is the wye-delta operation as in \cite{ar:delwye}, whose reversal is the delta-wye operation.

In general, any triminimal $2s$-regular graph $J$ can be transformed into a 
$(3,s)$-biregular bipartite graph $B$ by using delta-wye operations on each triangle in $J$. However, we can say more about the properties of $B$:

\begin{cor}\label{c:declaw}
Applying the wye-delta operation on all vertices in the first part of a $(3,s)$-biregular graph of girth greater than 6 creates a triminimal $2s$-regular graph.
\end{cor}

\begin{proof}
Triminimal 4-regular graphs are characterised this way in Corollary \ref{c:opt4}, noting that  the subdivision of a triangle in the 3-regular graph will produce a 6-cycle in the $(3,2)$-biregular graph. Additionally $(3,1)$-biregular graphs are simply the graphs $x K_{1,3}$ which have wye-delta transformation into the triminimal graphs $xK_3$ as required.

Suppose now that $s\geq 3$ and $J$ is a triminimal $2s$-regular graph with $n$ vertices and $ns$ edges.  
As in Corollary \ref{c:opt4}, $J$ has exactly $s$ triangles at each vertex and, using the delta-wye operation on each triangle in $J$, this gives $B$ which is a $(3,s)$-biregular graph.
A cycle of length 6 in $B$ would correspond to a triangle in $J$ beyond those 
guaranteed by the construction and cycles of length 4 would come from multiple edges in $J$, contrary to simplicity. 
\end{proof}

Without using $(3,s)$-biregular graphs it is also sometimes possible to create a triminimal $(2s+2)$-regular graph from a $2s$-regular triminimal graph $G$ by adding a 2-factor containing triangles from the complement of $G$.

\begin{figure}[h]
\begin{center}
\includegraphics[width=0.9\textwidth]{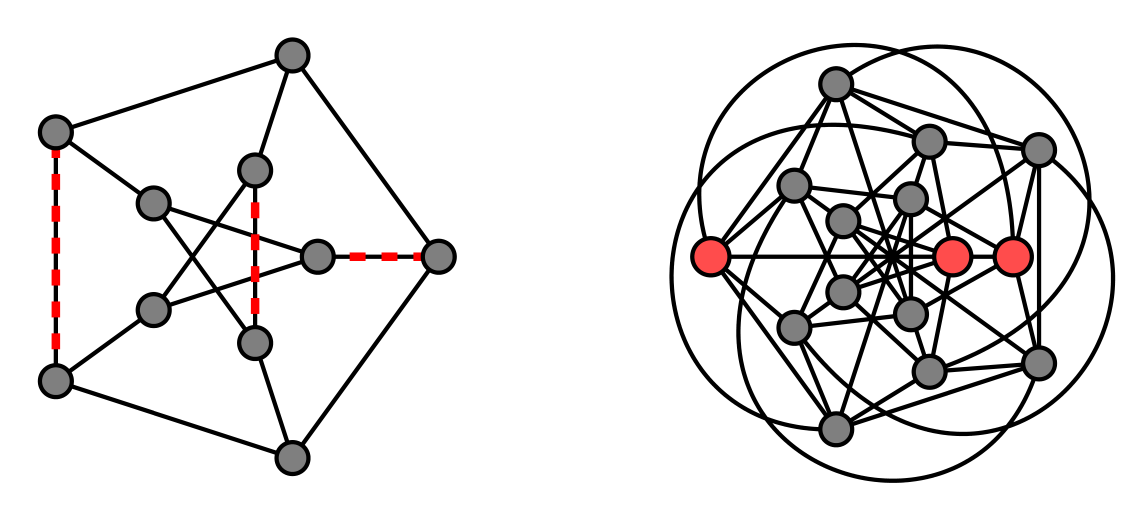}
\end{center}
\caption{Edges in Petersen at distance 3 giving a 6 regular graph with 15 triangles}\label{f:pete}
\end{figure}

For instance, as shown in red on the left of Figure \ref{f:pete}, in the Petersen graph there are five rotationally symmetric sets of three edges. If we take the line graph there are then five sets of three vertices distance 3 from each other. One set of three is shown in the figure on the right and the other sets are rotations of these. Together they can be used to make five new triangles without creating any other triangles including edges from the line graph. In this way we can create a 6-regular graph with $\frac{rn}{6}=\frac{15 \times 6}{6} = 15$ triangles. 
This graph also comes from applying Corollary \ref{c:declaw} to Tutte's 8-cage which was introduced in \cite{ar:T8}.

\subsection{Minimal but not triminimal graphs}

For the values of $n$ and $r$ which do not give integer values for $\lfloor\frac{r+1}{2}\rfloor\times \frac{n}{3}$, it is possible to produce graphs with only slightly more triangles than would be in a triminimal graph as is demonstrated here for $r=4$.

\begin{thm}
The 4-regular graphs with $n \equiv j ~({\rm mod}~ 3)$ vertices and the triangle property with the fewest triangles have $2\times \lfloor\frac{n}{3}\rfloor + 2$ triangles when $j\not = 0$.
\end{thm}

\begin{proof}
From Corollary \ref{c:trimin} we know when $r=4$ there are at least $\lceil\frac{2n}{3}\rceil$ triangles in any graph with the triangle property. Define $j:\equiv  n ~({\rm mod}~ 3)$ and consider the values of $j \not= 0$, since when $j=0$ we have the triminimal graphs themselves.

Firstly suppose $j=2$, so $n:=3k+2$ and $k=\lfloor\frac{n}{3}\rfloor$. Take a triminimal 4-regular graph with 
$3k$ vertices and $2k$ triangles and remove one triangle $T$ and add two vertices adjacent to each other and all three vertices of $T$, as in Figure \ref{f:nottri} (operation 2 in \cite{ar:PR}). This creates a 4-regular graph $H$ with $3k+2$ vertices and every edge in at least one triangle. In fact, all edges apart from the one between the two added vertices are in exactly one triangle. Thus there are $2k - 1 + 3 = 2k+2$ triangles and $\lceil\frac{2(3k+2)}{3}\rceil = 2k +2$ so $H$ has the triangle property and the fewest number of triangles possible.

\begin{figure}[h]
\begin{center}
\includegraphics[width=0.4\textwidth]{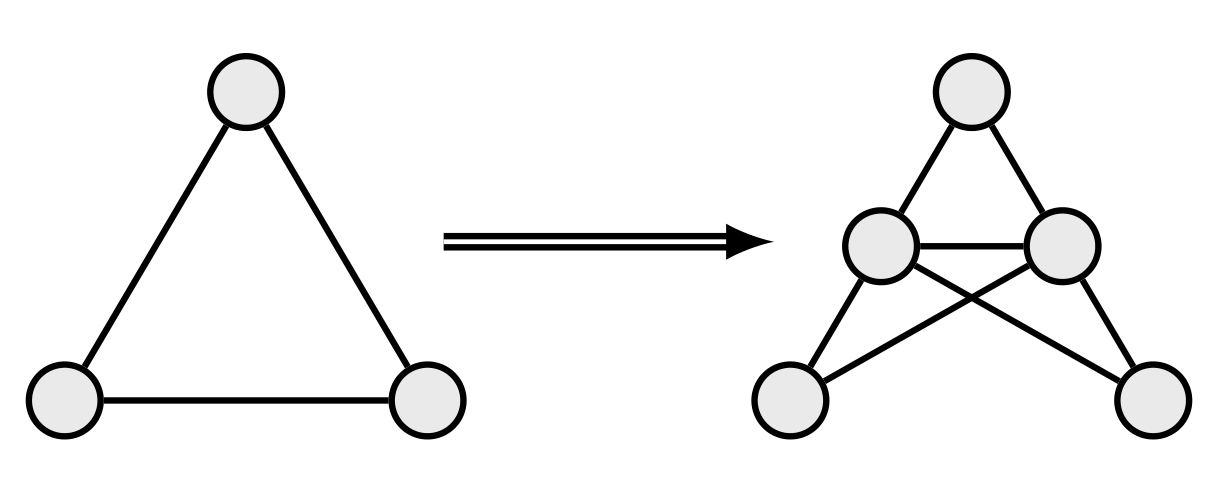}
\end{center}
\caption{Replacing a triangle by two new vertices}\label{f:nottri}
\end{figure}

We can proceed similarly when $n:=3k+1$, but this time it is necessary to start with a triminimal 4-regular graph with $3k-3$ vertices and $2k-2$ triangles; we then remove two triangles and add two pairs of vertices as in Figure \ref{f:nottri}. The resulting graph $J$ has $3k+1$ vertices and $(2k-2) -2 + 6 = 2k +2$ triangles. Note that $\lceil\frac{2(3k+1)}{3}\rceil = 2k +1$, but the only way this number of triangles could be achieved would be if there was exactly one vertex $v$ in three triangles, and that is not possible since at least one neighbour of $v$ will necessarily also be in three triangles. Hence $J$ is an $n$-vertex 4-regular graph with the smallest number of triangles 
\end{proof}

Simlar local operations can be defined for larger values of $r$, although most require several triangles beyond $\lfloor\frac{r+1}{2}\rfloor\times \frac{n}{3}$ and more cases are required. A particular feature making the case $r=4$ work well is that $r-1$ is a multiple of 3 and so we can remove $\frac{r-1}{3}$ triangles and add two vertices joined by an edge to give the new vertices degree $r$ as in Figure \ref{f:nottri}. Additionally, there also exist multigraphs with the triangle property and a low number of triangles, which further complicates matters, so further work needs to be done to fully clarify this area.

\section{Triminimal quintic graphs}\label{s:lgraph}

For 5-regular graphs, the situation is more complicated than in Section \ref{s:even}. The structure of the family of 5-regular graphs with the triangle property is investigated in an upcoming paper by the author but, unlike in \cite{ar:PR}, the base graphs can contain multiple edges and do not have the minimal number of triangles. 

\begin{thm}\label{t:triv}
Given a 5-regular graph with the triangle property, if it has $2n$ vertices then there are at least $2n$ triangles and also $n$ edges in at least two triangles.
\end{thm}

\begin{proof}
Suppose $v$ is a vertex in a 5-regular graph $G$ with $2n$ vertices and the triangle property. By Corollary \ref{c:trimin} $G$ contains at least $\frac{6}{2}\times \frac{2n}{3} =2n$ triangles and, by Theorem \ref{t:simple}, has no multiple edges. We can use Theorem \ref{t:mintri} to show that $v$ is incident with at least three triangles, and so $v$ is incident with at least one edge in two triangles or more. Since at all $2n$ vertices there is at least one such edge, there are at least $n$ edges of $G$ in at least two triangles.
\end{proof}

It turns out that there do exist graphs for which the inequality is tight in Theorem \ref{t:triv}, such as by the following construction of a quintic graph $L_n$, for $n\geq 7$:
\begin{enumerate}[label=(\alph*)]
  \item Create $n$ vertices labelled 0 to $n-1$ with edges joined in a cycle.
\item 
For all $k$ from 0 to $n-1$, join vertex $k$ to a new vertex labelled $n+k$.
\item For $k$ from 0 to $n-1$ add an edge between vertex $n+k$  and
\begin{enumerate}[label=(\roman*)]
\item vertex $n+((k+2)\mod n)$.
 \item vertex $((k-2)\mod n)$ 
\item vertex $((k+1)\mod n)$.
\end{enumerate}
\end{enumerate}

\begin{figure}[h]
\begin{center}
\includegraphics[width=0.5\textwidth]{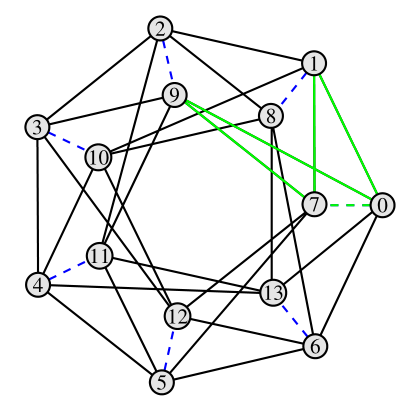}
\end{center}
\caption{The unique quintic graph $L_7$ containing 14 vertices and 14 triangles}\label{f:l7}
\end{figure}

$L_n$ is simple and 5-regular for $n\geq 5$; vertices 0 to $n-1$ are joined to two others in the cycle in (a) and one in (b) and two in (c), and the remaining $n$ vertices are joined to one vertex in (b) and four via (c). In $L_n$ (for $0\leq k \leq n-1$) the edges from vertex $k$ to vertex $n+k$, shown as the radial dashed spokes in Figure \ref{f:l7} for $n=7$, are in exactly two triangles, with vertex $(k+1) \mod n$ and vertex $n+ ((k+2)\mod n)$. The other two edges from vertex $k$ are to $(k-1)\mod n$ and $n+ ((k-1)\mod n)$; these vertices have an edge between them, so those edges are in a triangle. Similarly, vertex $n+k$ is adjacent to vertices $(k-2)\mod n$ and $n+ ((k-2)\mod n)$ and so the edges from vertex $n+k$ are part of that triangle. All of the different types of edges are therefore in a triangle and so $L_n$ has the triangle property; there are only three triangles at vertices $k$ and $n+k$ and since $L_n$ is formed symmetrically, that is true for all vertices of $L_n$.

Note that this implies that $L_n$ contains exactly $2n$ triangles too, by counting the triangles at each vertex and dividing by the number of vertices in each triangle.
Additionally these triangles partition $E(L_n)$ into $n$ edge-disjoint subgraphs isomorphic to the {\em diamond} ($K_{2,1,1}$). One of the seven edge-disjoint $\dtri$s $[0,1,7,9]$ is highlighed in Figure \ref{f:l7}; the edge 07 is in two triangles.

If $n\leq 6$ then extra triangles are formed in $L_n$ using edges between vertices $n$ to $2n-1$, and for $n\leq 4$ multiple edges are created by this construction, so the graphs created are not triminimal. $L_7$ is the smallest simple quintic triminimal graph, as shown by an computer search using nauty \cite{ar:nauty}.

As in Section \ref{s:even}, examples of triminimal $r$-regular graphs for larger odd $r$ can be constructed from a triminimal 5-regular graph $F$ with $2n$ vertices by carefully adding triangles between triples of vertices of $F$ so as not to create any triangles other than those we are choosing to add. Note that this is only possible, by Corollary \ref{c:trimin}, if $n$ is a multiple of 3 or $r$ is congruent to 2 mod 3.

For instance, if we take $L_{3j}$ for $j\geq 5$ and add edges between the triples of vertices between 0 and $3j-1$ based on their congruence mod $j$, and then similarly join triples of vertices between $3j$ and $6j-1$ based on their value mod $j$, we will only create those $2j$ extra triangles, and so we have created a family of triminimal 7-regular graphs.

\section{Construction of all triminimal quintic graphs}

\subsection{Diamonds and even-cycle decompositions}

\begin{figure}[ht]
\begin{center}
\includegraphics[width=1\textwidth]{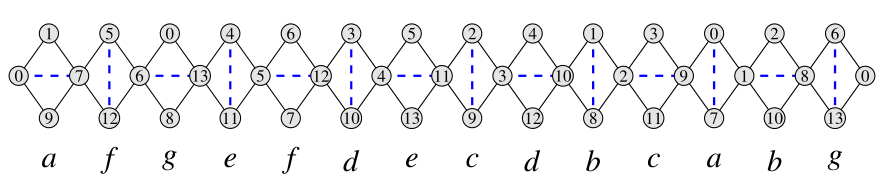}
\end{center}
\caption{Double cover of $E(L_7)$ by a sequence of diamonds}\label{f:l7tri}
\end{figure}

Any 5-regular graph with $2n$ vertices and $2n$ triangles must, as in Theorem \ref{t:triv}, be decomposable into $n$ edge-disjoint $\dtri$s.
In the figures in this section the edges that are in two triangles will again be coloured blue and dashed to highlight them, and they will be referred to as the {\em rotor} of the diamond. 
Moreover, we can cover all edges of $L_7$ twice using a cycle
of $\dtri$s, as shown in Figure \ref{f:l7tri}; if we label each $\dtri$ cyclically from $a:=\{0,1,7,9\}$, $b:=\{1,2,8,10\}$ to $g:=\{6,0,13,8\}$, we get the sequence indicated underneath the figure.

This gives, in general, {\em diamond sequences} for any triminimal quintic graph $G$; we pick any diamond $D_1$ and a vertex $v_1$ from the rotor of $D_1$ to start the sequence. Next we identify the other diamond that $v_1$ is in, say $D_2$, and $v_1$ must not be in the rotor of $D_2$. We can then uniquely identify $v_2$ as the other vertex not in the rotor of $D_2$ and find the other diamond $v_2$ is in, and so on, until the rotor of $D_1$ is reached. Note that in Figure \ref{f:l7tri} each diamond occurs once with the rotor vertical and once horizontal, and we will use the convention that the diamonds with the horizontal rotor are always located in the odd numbered positions in the sequences.
 
We can equivalently think about diamond sequences for $G$ by contracting each rotor edge and removing any multiple edges thus formed. This will give a 4-regular graph $X(G)$ with $\frac{|V(G)|}{2}$ vertices and the diamond sequences give rise to an {\em even-cycle decomposition} of $X(G)$; that is a partition of $E(X(G))$ into cycles of even length as in \cite{ar:ecd}. For $L_7$ we will have $X(L_7)$ as the complement of the 7-cycle, and the contraction of the rotors in Figure \ref{f:l7tri} gives the even-cycle $[0,5,6,4,5,3,4,2,3,1,2,0,1,6]$. The double cover given below the figure is exactly the same after substituting $a=0$, $b=1$,\ldots, $g=6$.

\begin{algorithm}\label{a:diam}
If we are given an even-cycle decomposition $C$ of a 4-regular graph with its vertices labelled by letters, we can follow the algorithm below to convert it into a 5-regular graph.
\begin{enumerate}
\item
Create vertex-disjoint copies of the diamond for each letter in $C$.
\item
Choose a cycle $Y$ of length $2j$ from $C$ and repeat steps 3 and 4 for 
$k:=1,\ldots, j$.
\item 
For the diamond in position $2k-1$ of $Y$ identify a rotor vertex from it with a non-rotor vertex in the diamond in position $2k$ of $Y$. 
\item
Now identify the other non-rotor vertex in diamond $2k$ with a rotor vertex in the diamond in position $(2k+1)$ of $Y$ (if $k=j$ then use the initial vertex in $Y$).
\item
Repeat step 2 for every other cycle in the even-cycle decomposition.
\end{enumerate}
\end{algorithm}

\begin{figure}[ht]
\begin{center}
\includegraphics[width=1\textwidth]{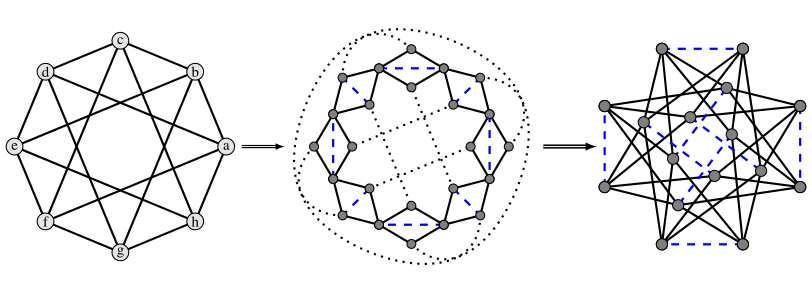}
\end{center}
\caption{Converting a 4-regular graph to a triminimal graph}\label{f:two8s}
\end{figure}

For example, given the decomposition $[[a, b, c, d, e, f, g, h], [b, g, d, a, f, c, h, e]]$, we can convert this particular pair of 8-cycles into the triminimal quintic graph shown in Figure
\ref{f:two8s}. We first create the diamonds and join them as per the first 8-cycle of the decomposition; the dotted lines in the second figure indicate the pairs of vertices to be identified to create the triminimal graph as per the second 8-cycle.

The even-cycle decomposition is given as a sequence of sequences of vertices.
In the example of Figure \ref{f:two8s}, $[g,d]$ is a subsequence of the second sequence, but $[a,c]$ is not since nowhere in either sequence is $c$ immediately preceded by $a$. 
Additionally, $[a,d]$ is not a subsequence since letter order is important but $[e,b]$ is a subsequence of the second sequence since we suppose each sequence is representing a cycle.

\begin{lem}\label{l:mint}
If G is a 4-regular simple graph with an even-cycle decomposition sat-
isfying the restriction that, for any vertices $x$, $y$ and $z$ in the even-cycle decomposition of the the 4-regular graph:
\begin{enumerate}[label=(\roman*)]
\item
Each vertex appears once in an odd numbered position and the other in an even position;
\item
Either $[x,y]$ or $[y,x]$ can appear as a subsequence, not both, and only once;
\item
If $[x,y]$ and $[y,z]$ both appear as subsequences then $[z,x]$ cannot appear as a subsequence unless all appear as a subsequence $[x,y,z,x]$;  
\end{enumerate}
then the application of Algorithm \ref{a:diam} yields a triminimal
5-regular graph. Conversely, any 5-regular triminimal graph can be
obtained in this fashion.
\end{lem}

\begin{proof}

Following Algorithm \ref{a:diam} will guarantee the creation of a 5-regular graph with the triangle property since every edge comes from a diamond and diamonds have the triangle property,
and the parity property (i) is necessary for 5-regularity.
We will show that neither of the properties (ii) or (iii) in the lemma can occur without creating edges in more than one triangle other than the rotors. This cannot occur in a triminimal graph and these edges are shown as dotted in Figure \ref{f:moretri}.

\begin{figure}[ht]
\begin{center}
\includegraphics[width=0.8\textwidth]{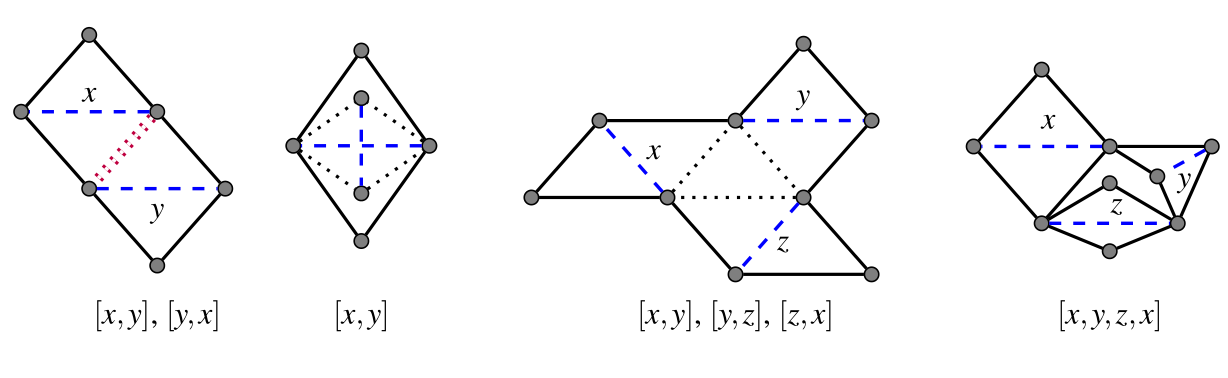}
\end{center}
\caption{Dotted edges are in more than one triangle so the graph cannot be triminimal}\label{f:moretri}
\end{figure}

In the left figure we have $[x,y]$ and $[y,x]$ appearing in two different sequences and note that the dotted edge is a double edge contributed from both diamonds and 
therefore in two extra triangles with the rotors of $x$ and $y$ so the resulting graph is not triminimal. 
If we had a sequence of length 2 such as $[x,y]$ then the second structure must exist and there are four edges in more than one triangle other than the rotors. Note that this would have to come from a double edge in the 4-regular graph.

Similarly, if $[x,y]$, $[y,z]$ and $[z,x]$ all appear in sequences but are not in the form $[x,y,z,x]$ then we can show there is a triangle formed by edges of the three diamonds as shown by the dotted edges in the third structure in Figure \ref{f:moretri}. This triangle is not part of one of the original $n$ diamonds, contradicting triminimality.

Within sequences, using property (i), we cannot have $[x,y,z]$ and $[z,x]$ as subsequences since $x$ and $z$ are the same parity positions in the former and the opposite parity in the latter. 
If $x$ is in an even position, say, then $[x,y]$ appearing means that $y$ is then in an odd position there, and so, if $[y,z]$ appears non-consecutively, $y$ must then be in an even position too, and then similarly for $z$ in $[z,x]$, and we get the dotted triangle as shown. If $x$, $y$ and $z$ are in odd positions instead we get the pattern shown as $[x,z]$, $[z,y]$ and $[y,x]$, similarly. 

Any longer sequences formed by the diamonds will not create any extra triangles, though, since the only way to get such triangles in the quintic graph is to have them in the 
even-cycle decomposition.
Note that for $L_n$ we do get $[x,y,z,x]$ appearing as a sequence, as shown in the fourth structure, and no extra triangles are formed, but they are consecutive diamonds as per the exception to the property. 
\end{proof}

\subsection{Removable edges}
One useful feature of any triminimal 5-regular graph $G$ is that all edges $e\in E(G)$ have the property that $G-e$ does not have every edge in a triangle. We shall call such an $e$ an {\em unremovable} edge and all other edges {\em removable}.

\begin{thm}\label{t:elig}
The only quintic graphs with the triangle property and with all edges as unremovable are the triminimal graphs.
\end{thm}

\begin{proof}
Suppose $G$ is a quintic graph with the triangle property and all of the edges of $G$ are unremovable. $G$ must be simple since an edge in a multiple edge is removable. 
Now we suppose $G$ is not triminimal, so has more triangles than vertices; therefore there must be more than $\frac{|V(G)|}{2}$ edges in $G$ in at least two triangles and hence there
exist two such edges $e_1$ and $e_2$ with a vertex $v$ in common. Because $G$ is simple and quintic and $e_1$ and $e_2$ are both in at least two triangles, either $e_1$ or $e_2$ are in a triangle with each other, or there is another edge (also in two triangles) incident with $v$ and in a triangle with one of $e_1$ or $e_2$.

\begin{enumerate}[label=\alph*)]
\item
Let us first suppose that $G$ does not contain a $K_4$.
Without loss of generality we can assume that $e_1$ and $e_2$ are two edges in two triangles that are also part of a triangle $T$ with each other. Moreover, the third edge of $T$ is also in at least two triangles since otherwise it is removable, and we therefore have the symmetrical situation shown in Figure \ref{f:blueadj} where $T$ is the triangle with the blue dashed edges.

\begin{figure}[ht]
\begin{center}
\includegraphics[width=0.4\textwidth]{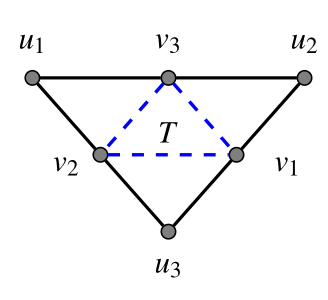}
\end{center}
\caption{No $K_4$ and adjacent edges in more than one triangle}\label{f:blueadj}
\end{figure}

 If $u_i = u_j$ (for some $i\not =j$), or $u_jv_j$ is an edge then we have $K_4$ as a subgraph, contrary to our supposition. We can thus assume that all 6 vertices are distinct and have no edges from any $v_j$ to these vertices other than those shown in Figure \ref{f:blueadj}.

All three vertices $v_1$, $v_2$ and $v_3$ in Figure \ref{f:blueadj} are degree 5 in $G$, and so need one more edge from them and that edge must be in a triangle.
If all are adjacent to the same vertex then there is a $K_4$, so we can assume that $v_3$, say, is adjacent to a vertex, $w_3$ and $v_1w_3$ and $v_2w_3$ are not edges of $G$. There must be an edge from $w_3$ to either $u_1$ or $u_2$ to make a triangle, and there is symmetry, so we can assume $w_3u_2$ is an edge, but this implies that $v_3u_2$ is in two triangles. Similarly we can now deduce that $u_2v_1$ is in two triangles too, as otherwise it is removable. However, this means that $v_1v_3$ is removable, a contradiction to all edges being unremovable.

\item
If there is a $K_4$ in $G$ then each edge of it is in two triangles in the $K_4$ alone. Moreover, each is removable unless it is part of another triangle. However, there are only two more edges from each vertex to other vertices in $G$, so at least one of the edges from each vertex in the $K_4$ is also in at least two triangles. This means that there must now be a subgraph isomorphic to $K_5 \backslash E(K_2)$, such that all edges in it are in at least two triangles. Now, as before, each of these edges will be removable unless they are part of a triangle with a new vertex, telling us that there must be a subgraph isomorphic to $K_6\backslash E(K_3)$. However, the edges between the vertices of degree 5 in this subgraph are removable.
\end{enumerate}
\end{proof}

\subsection*{Acknowledgements}
Many thanks to the referees for their very helpful suggestions which have greatly improved this paper.

\end{document}